\input amstex
\loadbold
\documentstyle{amsppt}
\magnification=\magstep1
\pagewidth{6.5truein}
\pageheight{9.0truein}

\hyphenation{semi-prime}


\def\xspec{X\operatorname{-spec}}

\def\gspec{G\operatorname{-spec}}
\def\gmax{G\operatorname{-max}}
\def\G{{\Cal G}}
\def\ggspec{{\Cal G}\operatorname{-spec}}
\def\ggmax{{\Cal G}\operatorname{-max}}

\def\ctil{\widetilde{c}}
\def\Rtil{\widetilde{R}}
\def\Atil{\widetilde{A}}
\def\phtil{\widetilde{\Phi}}

\def\ZZ{{\Bbb Z}}

\def\CC{{\Bbb C}}

\def\Hom{\operatorname{Hom}}
\def\kx{k^\times}
\def\chr{\operatorname{char}}

\def\spec{\operatorname{spec}}
\def\prim{\operatorname{prim}}
\def\max{\operatorname{max}}

\def\rad{\operatorname{rad}}

\def\O{{\Cal O}}
\def\bflam{{\boldsymbol \lambda}}

\def\spspec{S^{\perp}\text{-}\operatorname{spec}}
\def\spmax{S^{\perp}\text{-}\operatorname{max}}
\def\kgam{k\Gamma}
\def\sperp{S^{\perp}}

\def\swperp{S_w^{\perp}}
\def\swpspec{S^{\perp}_w\text{-}\operatorname{spec}}

\def\Oq{{\Cal O}_{\bold q}}

\def\bfq{{\bold q}}


\def\AbrHae{{\bf 1}}
\def\ArScTa{{\bf 2}}
\def\BroGoo{{\bf 3}}
\def\DeCKaPr{{\bf 4}}
\def\GolMic{{\bf 5}}
\def\GooLet{{\bf 6}}
\def\GooWar{{\bf 7}}
\def\HodLevone{{\bf 8}}
\def\HodLevtwo{{\bf 9}}
\def\HLT{{\bf 10}}
\def\Jos{{\bf 11}}
\def\Jbook{{\bf 12}}
\def\McCPet{{\bf 13}}
\def\Nor{{\bf 14}}
\def\Van{{\bf 15}}

\topmatter

\comment
\pretitle{
\hbox{}
\vskip -0.9truein
\hfill {\sevenrm August 26, 1998}
\vskip 0.9truein
         }
\endcomment
\title Quantum $n$-space as a quotient of classical $n$-space
 \endtitle

\rightheadtext{Quantum $n$-space}

\author K. R. Goodearl and E. S. Letzter \endauthor

\address Department of Mathematics, University of California, Santa
Barbara, CA 93106 \endaddress

\email goodearl\@math.ucsb.edu\endemail

\address Department of Mathematics, Texas A{\&}M University, College
Station, TX 77843\endaddress

\email letzter\@math.tamu.edu\endemail


\abstract The prime and primitive spectra of $\Oq(k^n)$, the
multiparameter quantized coordinate ring of affine $n$-space over an
algebraically closed field $k$, are shown to be topological quotients
of the corresponding classical spectra, $\spec \O(k^n)$ and $\max
\O(k^n)\approx k^n$, provided the multiplicative group generated by
the entries of $\bfq$ avoids $-1$.  \endabstract

\thanks The research of the first author was partially supported by
NSF grant DMS-9622876, and the research of the second author was
partially supported by NSF grant DMS-9623579. \endthanks

\endtopmatter
\document

\head Introduction \endhead

In the representation theory of a noncommutative ring $A$, the natural
analog of the maximal spectrum of a commutative ring is $\prim A$, the
set of primitive ideals (i.e., annihilators of irreducible
$A$-modules). Moreover, $\prim A$ can be equipped with the Zariski
topology (which in this setting is also called the Jacobson
topology). Thus, if $A$ is a quantization of a classical coordinate
ring $\O(V)$ over an algebraically closed field, one can view $\prim
A$ as a quantization of the variety $V$. The question then naturally
arises, how are $\prim A$ and $V$ related?  Some ``piecewise''
relations are known in many cases. For instance, if $G$ is a
connected, semisimple, complex algebraic group and $H$ is a maximal
torus of $G$, then various (generic) quantizations $A$ of $\O(G)$
exhibit the following properties: $H$ acts on $A$ via automorphisms;
there are only finitely many $H$-orbits in $\prim A$, and they are
locally closed; each $H$-orbit in $\prim A$ is homeomorphic to a
torus; and each $H$-orbit in $\prim A$ is a set-theoretic quotient of
a locally closed subset of $G$ that is stable under translation by $H$
(see \cite{\HodLevone, \HodLevtwo, \Jos, \Jbook, \HLT}). Similar
pictures have been observed to hold for quantized coordinate rings of
affine spaces (see \cite{\DeCKaPr, \Van, \BroGoo, \GooLet}). One is
therefore led to conjecture that in the above situations, $\prim A$ as
a whole is a topological quotient of $V$. Similarly, $\spec A$ (see
below) should be a topological quotient of $\spec \O(V)$. In the
``smallest'' cases, such as $A= \O_q(SL_2(k))$ and $A= \O_q(k^2)$,
these conjectures can easily be verified by direct calculation. In
this paper we establish these conjectures for multiparameter quantum
affine spaces, when $-1$ cannot be written as a product of the
defining parameters.

To provide some detail, first let $k$ be an algebraically closed
field. By a quantum $k$-affine $n$-space (i.e., a quantization of the
coordinate ring of $k^n$), we mean an algebra with generators
$x_1,x_2,\ldots,x_n$ subject only to the commutator relations $x_ix_j
= q_{ij}x_jx_i$, for a chosen set of nonzero scalars $q_{ij} \in
k$. (To avoid degeneracies, one assumes that $\bfq = (q_{ij})$ is a
multiplicatively antisymmetric $n\times n$ matrix.)  We denote this
algebra as $\Oq (k^n)$. Letting $H \colon= (\kx)^n$ be an algebraic
$n$-torus, it is not hard to verify that the natural action of $H$ on
the vector space spanned by $x_1,\ldots,x_n$ extends to an action by
$k$-algebra automorphisms on $\Oq (k^n)$. 

Next, set $A =\Oq (k^n)$. The $H$-action on $A$ described in the
previous paragraph extends naturally to actions on $\prim A$ and on
the larger set $\spec A$ of prime ideals of $A$. (Recall that a ring
is prime when every product of nonzero ideals is nonzero, and
that an ideal $I$ of $A$ is prime when $A/I$ is a prime ring. It is
important to note, for the prime ideals $I$ we are considering in this
paper, that $A/I$ may have zero-divisors.) We can impose the Zariski
topology on $\spec A$, and $\prim A$ then becomes a subspace under the
relative topology.

Our methods now require us to make the following technical assumption:
Either the subgroup of $\kx$ generated by the $q_{ij}$ does not
contain $-1$, or the characteristic of $k$ is $2$. (While our present
techniques do not appear to work in the absence of this hypothesis, we
do not know if our main results will remain valid without it.) Now let
$R\colon = \O(k^n)$ be a commutative polynomial ring in $n$ variables,
also equipped with the natural action of $H$ by automorphisms, and
continue to let $A = \Oq (k^n)$. In the main result (4.11),
we construct $H$-equivariant topological quotient maps
$$ \spec R \rightarrow \spec A \quad \text{and} \quad k^n \approx \max
R \rightarrow \prim A .$$
Moreover, a description of the fibers of the second map is provided.
Roughly speaking, we show that $\prim A$ is equal to $k^n$ modulo the
actions of a compatible system of groups acting separately on each
$H$-orbit, and that $\spec A$ can be obtained from $\spec R$ in an
analogous fashion. These quotient maps are not uniquely determined by
the $q_{ij}$, but depend only on a choice of square roots in the
construction of a specific cocycle $c$; see (3.12).

Our analysis involves two separate steps: 
\roster 
\item"(I)" The
spaces $\spec A$ and $\prim A$ are homeomorphic to certain sets of
ideals of $R$, denoted $\ggspec R$ and $\ggmax R$, which are equipped
with Zariski topologies.
\item"(II)" We obtain topological quotient maps $\spec R \rightarrow
\ggspec R$ and $\max R \rightarrow \ggmax R$, by sending each prime
(respectively, maximal) ideal $P$ of $R$ to the largest member of
$\ggspec R$ contained within $P$.
\endroster
Step II works very generally, and so we begin with that. Section 1
contains an axiomatic treatment that produces topological quotients of
prime and primitive spectra in a noetherian ring $R$, using certain
families of ideals of $R$. (When we say that a ring is noetherian, we
mean left and right noetherian.) In Section 2, we restrict to the case
of a commutative polynomial ring $R$ over $k$ and a family of ideals
stable under a suitably compatible family of groups of automorphisms
of $R$. This last setup yields the topological quotients of $\spec R$
and $\max R$ that will appear in our main theorem. Step I is developed
in Sections 3 and 4. Section 3 is devoted to the case of quantum tori,
which yields the necessary information about individual ``strata'' of
prime and primitive ideals in $A$; in Section 4, we glue these results
together over the whole prime and primitive spectra of $A$, and obtain
the main theorem. Section 5 contains explicit calculations of the
resulting topological quotient maps $k^n \rightarrow \prim A$ for the
case where $A$ is a single parameter quantum affine space, and Section
6 applies the main theorem to twists of more general graded
commutative algebras.
\medskip

We remind the reader that a map $\phi : X\rightarrow Y$ between
topological spaces is a {\it topological quotient map\/} provided
$\phi$ is surjective and the topology on $Y$ coincides with the
quotient topology induced by $\phi$, that is, a subset $W\subseteq Y$
is closed in $Y$ precisely when $\phi^{-1}(W)$ is closed in $X$. When
such a map exists, $Y$ is completely determined (up to homeomorphism)
by $X$ and the fibers of $\phi$.

\head 1. Quotients of prime spectra \endhead

We begin the work of this paper by specifying sufficient conditions
under which a collection of ideals in a noetherian ring $R$, equipped with a
Zariski topology, will form a topological quotient of the prime spectrum of
$R$. These conditions are modelled on standard properties of the collection of
those ideals of $R$ which are prime relative to the ideals invariant under
some set of operators, such as the collection of $G$-prime ideals where $G$ is
a group of automorphisms of $R$, or the collection of $\Delta$-prime ideals
where $\Delta$ is a set of derivations on $R$. Conditions for a subcollection
to be a topological quotient of the primitive or the maximal spectrum of $R$
are also given. At the end of the section, we specialize to the case of
$G$-prime ideals, this being the case we shall need for application to quantum
tori.

\subhead 1.1 \endsubhead Let $R$ be a noetherian ring, and let $\xspec
R$ be a nonempty set of ideals of $R$. The elements of $\xspec R$ will
be referred to as {\it $X$-prime\/} ideals, and intersections of
collections of
$X$-prime ideals will be termed {\it $X$-semiprime\/} ideals. We will
further assume, throughout this section, that

(a) For each prime ideal $P$ of $R$ there exists a (unique) $X$-prime
ideal, denoted $(P : X)$, that is maximum among $X$-semiprime ideals
contained in $P$. In other words, $(P:X)$ is an $X$-prime ideal contained in
$P$ and containing all $X$-semiprime ideals which are contained in $P$.

(b) If $P$ is a prime ideal of $R$ minimal over an $X$-prime ideal
$Q$, then $Q = (P : X)$. In particular, for each $X$-prime ideal $Q$
there exists a prime ideal $P$ such that $Q = (P : X)$.

\definition{1.2} If $I$ is an ideal of $R$ then $V(I)$ will denote the
set of prime ideals containing $I$. The set of prime ideals of $R$
will be denoted $\spec R$ and will be equipped with the standard
Zariski topology: The closed subsets are precisely those of the form
$V(I)$, for ideals $I$ of $R$. The sets of (left) primitive and
maximal ideals of $R$ will be denoted $\prim R$ and $\max R$
respectively, and each of these sets will be given the relative
topology from $\spec R$. \enddefinition

\proclaim{1.3 Lemma} Let $Q$ be an $X$-prime ideal in $R$, and let $I$
and $J$ be $X$-semiprime ideals of $R$ such that $IJ \subseteq
Q$. Then $I$ or $J$ is contained in $Q$. \endproclaim

\demo{Proof} Choose a prime ideal $P$ such that $Q = (P : X)$. Then
$IJ \subseteq P$, and so one of $I$ or $J$ is contained in $P$, say
$I\subseteq P$. But $Q$ is the unique maximum $X$-semiprime ideal contained in
$P$, and so $I \subseteq Q$. \qed\enddemo

\definition{1.4} If $I$ is an ideal of $R$, define $V_X(I)$ to be the
set of $X$-prime ideals of $R$ containing $I$; note that $V_X(I) =
V_X(J)$, where $J$ is the intersection of the $X$-prime ideals of $R$
containing $I$ (with the convention that the intersection of an empty
collection of ideals is equal to $R$ itself). By (1.3), if
$J_1,\ldots,J_t$ are $X$-semiprime ideals of $R$, then $V_X(J_1)\cup
\cdots \cup V_X(J_t) = V_X(J_1\cap \cdots \cap J_t)$.  Therefore, the
standard Zariski topology generalizes to a topology on $\xspec R$: The
closed sets are those subsets of the form $V_X(I)$, for ideals $I$
of $R$.  \enddefinition

\definition{1.5} By our assumptions in (1.1), there is a surjection
$$\pi: \spec R \rightarrow \xspec R$$
sending each prime ideal $P$ to $(P : X)$. We first record two trivial
observations concerning the behavior of $\pi$ with respect to closed sets.
\enddefinition

\proclaim{Lemma} Let $J$ be an $X$-semiprime ideal of $R$. 

{\rm (a)} $\pi (V(J))=V_X(J)$.

{\rm (b)} $\pi ^{-1}(V_X(J)) = V(J)$. \qed \endproclaim

\definition{1.6} A subset $U$ of $\spec R$ will be termed {\it
$\pi$-stable\/} provided $U$ is a union of fibers of $\pi$,
or equivalently, provided $U= \pi^{-1}(\pi(U))$. \enddefinition

\proclaim{Lemma} Let $U$ be a closed $\pi$-stable subset of $\spec
R$. Then there exists an $X$-semiprime ideal $J$ of $R$ such that $U =
V(J)$. \endproclaim

\demo{Proof} Choose an ideal $I$ of $R$ such that $U = V(I)$, and let
$P_1,\ldots,P_t$ be the prime ideals minimal over $I$. Then
$P_1,\ldots,P_t$ are the minimal elements of $U$. For $1 \leq i \leq
t$, set $Q_i = (P_i : X )$, and set $J = Q_1 \cap \cdots \cap
Q_t$. Observe that $U \subseteq V(J)$. Next, let $P$ be a prime ideal
minimal over $J$. Note that $P$ is minimal over $Q_i$, for some $1
\leq i \leq t$, and so $(P : X) =Q_i = (P_i : X)$. Therefore, since
$U$ is $\pi$-stable, $P \in U$. Consequently, $V(J) \subseteq U$, and
hence $U = V(J)$. \qed\enddemo

\proclaim{1.7 Proposition} {\rm (a)} $\pi$ is continuous.

{\rm (b)} $\pi$ maps $\pi$-stable closed subsets of $\spec R$ to closed
subsets of $\xspec R$.

{\rm (c)} $\pi$ is a topological quotient map. \endproclaim

\demo{Proof} (a) Let $V$ be a closed subset of $\xspec R$. As seen in
(1.4), there exists an $X$-semiprime ideal $J$ in $R$ such that $V =
V_X(J)$. By (1.5b), $\pi ^{-1}(V) = V(J)$.

(b) Let $U$ be a $\pi$-stable closed subset of $\spec R$. It follows
from (1.6) that there exists an $X$-semiprime ideal $J$ in $R$ such
that $U = V(J)$. By (1.5a), $\pi (U) = V_X(J)$.

(c) Assume that $V$ is a subset of $\xspec R$ whose inverse image $U =
\pi^{-1}(V)$ is closed in $\spec R$. It only remains to show that $V$ is
closed. However, $U$ is clearly $\pi$-stable, and so $V$ is closed by
(b).  \qed\enddemo

\definition{1.8} Let us say that a subset $U$ of $\prim R$ (or of
$\max R$) is {\it relatively $\pi$-stable\/} provided $U= V\cap \prim
R$ (or $U= V\cap \max R$) for some $\pi$-stable subset $V$ of $\spec
R$.  \enddefinition

\proclaim{Proposition} {\rm (a)} Suppose that $\bigcap U$ is
$X$-semiprime for every relatively $\pi$-stable subset $U$ of $\prim
R$. Then the restriction of $\pi$ to the surjection $$\prim R
\longrightarrow \pi(\prim R)$$ is a topological quotient map (with
respect to the relative topologies).

{\rm (b)} Similarly, if the intersection of any relatively
$\pi$-stable subset of $\max R$ is $X$-semi\-prime, then $\pi$
restricts to a topological quotient map $$\max R \longrightarrow
\pi(\max R).$$ \endproclaim

\demo{Proof} (a) Let $\pi_{\text{prim}} : \prim R \rightarrow
\pi(\prim R)$ denote the restriction of $\pi$. It is immediate from
the continuity of $\pi$ that $\pi_{\text{prim}}$ is continuous. Now
let $W$ be a subset of $\pi(\prim R)$ whose inverse image $U=
\pi_{\text{prim}}^{-1}(W)$ is closed in $\prim R$. Since $U$ is
relatively $\pi$-stable, its intersection is by hypothesis an
$X$-semiprime ideal $J$.  Moreover, $U= V(J)\cap \prim R$ because $U$
is closed in $\prim R$. Then $W= \pi(U)\subseteq \pi(V(J))=
V_X(J)$. If $I\in V_X(J)\cap \pi(\prim R)$, there exists $P\in \prim
R$ such that $\pi(P)=I\supseteq J$. Since $P\supseteq \pi(P)$, we have
$P\in V(J)\cap \prim R= U$, and so $I\in W$. Thus $W= V_X(J)\cap
\pi(\prim R)$, a relatively closed subset of $\pi(\prim R)$.

(b) Use the proof above, with ``prim'' replaced by ``max'' everywhere.
\qed\enddemo

The preceding observations apply easily to group actions, as follows.

\definition{1.9} For the remainder of this section, let $G$ be a group
acting by automorphisms on $R$. An ideal $I$ of $R$, such that $g(I) =
I$ for all $g \in G$, is termed a {\it $G$-ideal\/}. Note that the set
of $G$-ideals of $R$ is closed under products, sums, and
intersections. A (set-theoretically) proper $G$-ideal $Q$ of $R$
is {\it $G$-prime\/} provided $IJ \not\subseteq Q$ for all $G$-ideals
$I$ and $J$ not contained in $Q$. An intersection of $G$-prime ideals
will be termed {\it $G$-semiprime\/}. For any ideal $I$ of $R$, set
$(I : G) = \bigcap _{g \in G} g (I)$. It is easy to see that $(I : G)$
is the unique maximum $G$-ideal contained in $I$. \enddefinition

\proclaim{1.10 Lemma} {\rm (cf\. \cite{\GolMic, Remarks $4^*$,
$5^*$})} {\rm (a)} Let $P$ be a prime ideal of $R$. Then $(P : G)$ is
a $G$-prime ideal.

{\rm (b)} Let $Q$ be a $G$-prime ideal of $R$, and let $P$ be a prime
ideal of $R$ minimal over $Q$. Then $Q = (P : G)$, and every prime
ideal minimal over $Q$ is in the $G$-orbit of $P$. In particular, $Q$
is semiprime.
\endproclaim

\demo{Proof} Part (a) follows from an argument mimicking the proof of
(1.3).

(b) Let $P_1,\ldots,P_t$ be the prime ideals of $R$ minimal over $Q$,
and set $N_i = (P_i : G)$ for each $1 \leq i \leq t$. Note that
$N_1,\ldots, N_t$ are $G$-ideals, all containing $Q$, and that some
product of them is contained within $Q$. Therefore, $Q = N_i$ for some
$i$, because $Q$ is $G$-prime.  Now consider an arbitrary prime ideal
$P_j$ minimal over $Q$. Since $P_j \supseteq Q=N_i$, we see that
$P_j$ contains a product of prime ideals in the $G$-orbit of $P_i$,
and so $P_j \supseteq g(P_i) \supseteq Q$ for some $g \in
G$. Therefore, $P_j = g(P_i)$, and $Q = (P_i : G ) = (P_j : G)$. Part
(b) follows.  \qed \enddemo

\definition{1.11} Let $\gspec R$ denote the set of $G$-prime ideals of
$R$. It follows from (1.9) and (1.10) that $\gspec R$ satisfies the
axioms for $\xspec R$ specified in (1.1). In particular, if $V_G(I)$
denotes the set of $G$-prime ideals of $R$ that contain a given ideal
$I$ of $R$, then by (1.4) there is a Zariski topology on $\gspec R$:
The closed sets are those of the form $V_G(I)$, for ideals $I$ of
$R$. With respect to this topology, we can use (1.7c) to deduce that
the assignment
$$P \mapsto (P : G)$$
is a topological quotient map from $\spec R$ onto $\gspec R$.

Let $\gmax R$ denote the set of maximal proper $G$-ideals. Observe
that $\gmax R$ is a subset of $\gspec R$ and that $\gmax R$ is
comprised of the maximal members of $\gspec R$. Equip $\gmax R$ with
the relative topology. While each member of $\gmax R$ must be equal to
$(M:G)$ for some $M\in \max R$, the converse fails in general. For
example, let $R$ be a polynomial ring $k[x]$ over an infinite field
$k$, and let $G=\kx$ act on $R$ by the rule $g.f(x)= f(gx)$. Then
$(\langle x-\alpha \rangle :G)=0$ for all nonzero $\alpha \in k$,
whereas $\gmax R$ consists of the single maximal ideal $\langle
x\rangle$.

In the present setting, the hypothesis of (1.8b) is easily verified, as
follows. Suppose $U$ is a relatively $\pi$-stable subset of $\max R$, and set
$J= \bigcap U$. Since $\max R$ is stable under the action of $G$ (within
$\spec R$), so is $U$, and hence $J= \bigcap_{M\in U} (M:G)$. Thus by (1.10a),
$J$ is $G$-semiprime, as desired. We therefore conclude from (1.8b) that the
assignment
$$M\mapsto (M:G)$$
is a topological quotient map from $\max R$ onto $\pi(\max R)$.
\enddefinition

\definition{1.12} Suppose that $G$ has no proper subgroups of finite
index. Then all finite $G$-orbits of ideals of $R$ are singletons, and
so all $G$-prime ideals of $R$ are prime, by (1.10b). We see, in this
case, that $\gspec R$ is a subset of $\spec R$ and that the topology
on $\gspec R$ described in (1.11) is the relative topology.
\enddefinition

\definition{1.13} We conclude this section by considering a particular
case of the group action setting that will be needed later.

Let $k$ be an algebraically closed field, $R= k[y_1^{\pm1} ,\dots,
y_n^{\pm1}]$ a Laurent polynomial ring over $k$, and $H= (\kx)^n$ the
standard algebraic $k$-torus of rank $n$.  Let $H$ act on
$R$ by $k$-algebra automorphisms in the natural manner; namely,
$$(h_1,\ldots,h_n).f(y_1,\ldots,y_n) = f(h_1y_1,\ldots,h_ny_n),$$
for $(h_1,\ldots,h_n) \in H$ and $f(y_1,\ldots,y_n) \in R$. Since $k$ is
algebraically closed, the induced action of $H$ on $\max R$ is transitive.
\enddefinition

\proclaim{Proposition} Let $G$ be a subgroup of $H$.

{\rm (a)} The set $\{ (M:G) \mid M\in \max R\}$ coincides with $\gmax R$.

{\rm (b)} If $G$ is a closed subgroup of $H$, the fibers of the map
$M\mapsto (M:G)$, from $\max R$ to $\gmax R$, are precisely the
$G$-orbits in $\max R$. \endproclaim

\demo{Proof} (a) Given $Q\in \gmax R$, choose a maximal ideal $N\supseteq
Q$. Then $(N:G) \supseteq Q$, and so $Q= (N:G)$ by the maximality of $Q$.

Conversely, consider $Q= (M:G)$ where $M\in \max R$. Choose $Q'\in
\gmax R$ containing $Q$. By the preceding paragraph, there exists a
maximal ideal $M'$ in $\max R$ for which $Q' = (M':G)$. Since $H$ acts
transitively on $\max R$, there exists an $h \in H$ such that $M =
h(M')$. Since $H$ is abelian, $Q = (M : G) = (h(M') : G) = h(Q')$. 
Hence $Q' \supseteq h(Q')$, and so $Q'= h(Q')= Q$. Thus $Q\in \gmax R$.

(b) Let $M,M'\in \max R$. If $M$ and $M'$ are in the same $G$-orbit,
then $(M:G) = (M':G)$. Conversely, assume that $(M:G) = (M':G)$.  On
intersecting with the ring of $G$-invariants $R^G$, we thus obtain
$$M\cap R^G = (M:G)\cap R^G = (M': G)\cap R^G = M'\cap R^G.$$

Next, we may view $R$ as the coordinate algebra of $H$ and identify
$H$ with $\max R$. By \cite{\Nor, Theorem 48, p\. 220}, for example,
there is an identification of $H/G$ with $\max R^G$ that produces the
following commutative diagram:
$$\CD H @>\operatorname{quotient}>> H/G \\ @| @|\\ \max R
@>\operatorname{restriction}>> \max R^G \endCD$$
Letting $H$ act on itself and on $H/G$ by right translation (i.e., $g \in
H$ acts on $h \in H$ by $g{.}h= hg^{-1}$), the diagram above is
$H$-equivariant.

Since the cosets of $G$ in $H$ are the $G$-orbits with
respect to right translation, it follows that the fibers of the
restriction map $\max R \rightarrow \max R^G$ are the
$G$-orbits in $\max R$. Therefore $M$ and $M'$ lie in the
same $G$-orbit, as desired. \qed\enddemo

\head 2. Quotients of affine space \endhead

Throughout this section, $k$ denotes a field, and $R$ denotes a commutative
polynomial ring $k [y_1,\ldots,y_n]$. We study a ``piecewise'' action
on $R$ by a compatible system of groups, together with a collection of ideals
which we show satisfies the axioms of the previous section. The resulting
quotients of the prime and maximal spectra of
$R$ will play a crucial role in our main theorem. Namely, as we shall prove in
Section 4, the quotients of $\spec R$ and $\max R$ with respect to a suitable
system of groups turn out to be homeomorphic to the prime and primitive
spectra, respectively, of a quantum affine $n$-space.

The setup we develop here amounts to patching together finitely many quotients
by group actions, as follows. We first partition $\spec R$ into the $2^n$
locally closed sets determined by which subsets of $\{y_1,\dots,y_n\}$ are
contained in given prime ideals. These subsets are homeomorphic, via
localization, to the prime spectra of the Laurent polynomial rings obtained by
factoring some of the $y_i$ out of $R$ and inverting the remainder. Our
``piecewise group action'', finally, amounts to compatible choices of groups
acting as automorphisms on the above localizations.

\definition{2.1} Set $H = (\kx)^n$, the algebraic $k$-torus of rank
$n$, and equip $H$ with the Zariski topology. As in (1.13), let $H$
act on $R$ by the automorphisms
$$(h_1,\ldots,h_n).f(y_1,\ldots,y_n) = f(h_1y_1,\ldots,h_ny_n).$$
We will also consider the induced actions of $H$ on $\spec R$ and $\max
R$. Next, let $W$ denote the set of subsets of $\{1,\ldots,n\}$. For
each $w \in W$, let $\spec _w R$, equipped with the relative topology,
denote the set of prime ideals $P$ of $R$ such that
$$P \cap \{ y_1,\ldots,y_n\} = \{ y_i \mid i \in w \} .$$
Note that $\spec R = \bigsqcup _{w \in W} \spec _w R $.  Also note
that $\spec_w R$ is stable under $H$, and that there is an obvious
$H$-equivariant homeomorphism $$\spec _w R @>{\ \cong\ }>> \spec
k[y_i^{\pm 1} \mid i \not\in w].$$ \enddefinition

\definition{2.2} Let $\G = \{ G_w \mid w \in W\}$ be a family of
subgroups of $H$, indexed by $W$, subject to the following
compatibility hypothesis: $(P : G_v) \subseteq (P : G_w)$ for all $v
\subseteq w$ and every $P \in \spec_w R$. We do not exclude the
possibility that the $G_w$, for $w \in W$, are identically the
same. Set
$$\ggspec _wR \; = \; \left\{ (P : G_w) \mid P \in \spec _w R
\right\},$$
for $w\in W$. Note, since the $y_i$ are $G_w$-eigenvectors, that
$Q\cap \{y_1,\dots,y_n\}= \{y_i \mid i\in w\}$ for all $Q\in \ggspec_w
R$.  Consequently, the sets $\ggspec_w R$ for $w\in W$ are pairwise
disjoint, and we set
$$\ggspec R = \bigsqcup _{w \in W} \ggspec _wR .$$
The members of $\ggspec R$ will be termed {\it $\G$-prime\/} ideals, the
intersections of collections of $\G$-prime ideals will be called {\it
$\G$-semiprime\/} ideals, and the set of $\G$-prime ideals containing an
ideal $I$ of $R$ will be denoted $V_\G(I)$. Note from (1.10b) that all
$\G$-semiprime ideals are semiprime.  \enddefinition

\proclaim{2.3 Lemma} $\ggspec R$ satisfies the conditions in {\rm
(1.1)}, with $(P : \G) = (P : G_w)$ for $P \in \spec _w
R$. \endproclaim

\demo{Proof} To verify (1.1a), let $P \in \spec _w R$, let $Q = (P :
G_w)$, and let $Q'$ be an arbitrary $\G$-semiprime ideal of $R$
contained in $P$. We will show that $Q' \subseteq Q$. To start, write
$$Q' = \bigcap _{v\in W} N_v ,$$
where each $N_v$ is an intersection of ideals from $\ggspec _v R$, and
where each $N_v$ contains $Q'$. Since $Q'\subseteq P$, some
$N_v\subseteq P$. Then since $y_i\in N_v$ for all $i\in v$, we must
have $v\subseteq w$. Note that $N_v$ is $G_v$-stable and so is
contained in $(P : G_v)$. By the compatibility hypothesis, $N_v
\subseteq Q$, whence $Q'\subseteq Q$, and (1.1a) is established.

To check (1.1b), let $Q \in \ggspec _w R$ and let $P$ be a
prime ideal of $R$ minimal over $Q$. Observe first that $P \in \spec
_v R$ for some $v \supseteq w$. Next, it follows from (1.10a) that $Q$ is
$G_w$-prime, and then from (1.10b) that $Q = (P : G_w)$. Moreover, since
each of the variables $y_i$ is a $G_w$-eigenvector, it now follows that $v=w$,
whence $Q = (P : \G)$ and  (1.1b) is satisfied. \qed\enddemo

\definition{2.4} Following (1.4), we equip $\ggspec R$ with the
Zariski topology where the closed subsets are the $V_\G(I)$,
for ideals $I$ of $R$. In case none of the $G_w$ has proper subgroups of
finite index, it follows from (1.12) that $\ggspec R\subseteq \spec R$,
and then the topology on $\ggspec R$ coincides with the relative topology.

Note, for each $w \in W$, that $\spec _wR$ is a union of
$H$-orbits. Therefore, since $H$ is abelian, we see that $h(P : \G) =
(h(P) : \G)$ for all $h \in H$ and $P \in \spec R$. From (1.7) we now
obtain: \enddefinition

\proclaim{Proposition} The assignment $P \mapsto (P : \G)$ produces an
$H$-equivariant topological quotient map from $\spec R$ onto $\ggspec
R$. \qed\endproclaim

\definition{2.5} For $w \in W$, set $\max _w R = \max R \cap \spec _w
R$.  Note that $\max R = \bigsqcup _{w \in W} \max _w R$.  Moreover, $\max
_w R$ is equal to the set of maximal members of $\spec _wR$.  Next, set
$$\align \ggmax _wR \; &= \; \left\{ (M : G_w) \mid M \in \max _w R \right\}
\qquad (\text{for\ } w\in W);\\
\ggmax R \; &= \; \bigsqcup _{w \in W} \ggmax _wR = \{(M:\G) \mid M\in
\max R\}.
\endalign
$$
Note that the definition of $\ggmax_w R$ above and the
definition of $G_w\operatorname{-max} R$ in (1.11) follow different
patterns. However, we shall see in (2.6) that when $k$ is algebraically
closed,
$\ggmax_w R$ coincides with $\bigl( G_w\operatorname{-max} R \bigr) \cap
\bigl(
\ggspec_w R
\bigr)$.
\enddefinition

\proclaim{Proposition} The assignment $M \mapsto (M : \G)$
produces an $H$-equivariant topological quotient map from $\max R$
onto $\ggmax R$. \endproclaim

\demo{Proof} By (1.8b), it suffices to show that if $U$ is any relatively
$\pi$-stable subset of $\max R$ and $J= \bigcap U$, then $J$ is $\G$-semiprime.

To start, set $U_w= U\cap \max_w R$ for $w\in W$, and note that $U_w$ is
stable under $G_w$ since $(g(M):\G)= (M:\G)$ for all $M\in \max_w R$ and $g\in
G_w$. Hence,
$$\bigcap U_w= \bigcap_{M\in U_w} (M:G_w)= \bigcap_{M\in U_w} (M:\G),$$
which is $\G$-semiprime. Therefore $J= \bigcap_{w\in W} \bigl( \bigcap U_w
\bigr)$ is a $\G$-semiprime ideal, as required.
\qed
\enddemo

\proclaim{2.6 Proposition} Assume that $k$ is algebraically closed. For $w
\in W$, the set $\ggmax _wR$ coincides with the set of maximal members of
$\ggspec _wR$.
\endproclaim

\demo{Proof} Set $m = n - \vert w \vert$ and $R_w= \bigl( R/\langle
y_i\mid i\in w\rangle \bigr) [y_j^{-1}\mid j\notin w]$. Note that
$R_w$ is a Laurent polynomial ring in the (images of the) $m$
indeterminates $y_j$, for $j\notin w$. We thereby obtain an action of
the $m$-torus $(\kx)^m$ on $R_w$, following (1.13). Note that the
actions on $R_w$ by $G_w$ and $H$, induced from their actions on $R$,
factor through the $(\kx)^m$-action. Moreover, localization provides
$H$-equivariant bijections $$\max_w R \rightarrow \max R_w
\qquad\text{and}\qquad \ggspec_w R \rightarrow G_w\operatorname{-spec}
R_w.$$ The proposition therefore follows from (1.13a). \qed\enddemo

\definition{2.7} We close this section by considering the fibers of
the map in (2.5) more closely. First note that each $\max_w R$ is a
disjoint union of fibers. Next, on $\max_w R$ the map is given by the
rule $M\mapsto (M:G_w)$, and so the fibers within $\max_w R$ are
unions of $G_w$-orbits. Our interest is in conditions under which the
fibers within $\max_w R$ coincide with the $G_w$-orbits.

The maximal ideals in $\max_w R$ have the form
$$\langle y_i\mid i\in w\rangle +\langle f_1,\dots,f_t \rangle$$
where $f_1,\dots,f_t \in k[y_j \mid j\notin w]$. These ideals are all fixed
by the group
$$H_w= \{ (h_1,\dots,h_n)\in H \mid h_j=1 \text{\ for all\ } j\notin w\};$$
thus $(M:G_w)= (M:G_wH_w)$ for all $M\in \max_w R$. Hence, there is no loss of
generality in assuming that $G_w\supseteq H_w$.
\enddefinition

\proclaim{Proposition} Fix $w \in W$, and assume that $k$ is
algebraically closed. Assume that $G_w\supseteq H_w$ and that
$G_w/H_w$ is a closed subgroup of $H/H_w$. Then the fibers of the map
$M\mapsto (M:\G)$, from $\max_w R$ to $\ggmax_w R$, are precisely the
$G_w$-orbits in $\max_w R$. \endproclaim

\demo{Proof} We shift everything to the localization $R_w$, as in the
proof of (2.6), and we recall the $(\kx)^m$-action on $R_w$ described
therein, where $m = n - \vert w \vert$.  Now observe that the induced
action of $H/H_w$ on $R_w$ is identical to the $(\kx)^m$-action on
$R_w$. Therefore, the proposition follows from (1.13b).  \qed\enddemo

\head 3. Quantum tori \endhead

As in \cite{\GooLet}, much of our analysis of quantum affine spaces can be
reduced, via localization, to quantum tori. This portion of our analysis is
carried out in the present section. It leads, in particular, to a version of
our main theorem for an arbitrary multiparameter quantum torus $\Oq((\kx)^n)$
over an algebraically closed field $k$, which establishes that $\spec
\Oq((\kx)^n)$ and $\prim \Oq((\kx)^n)$ can be presented as topological
quotients of $\spec \O((\kx)^n)$ and $\max \O((\kx)^n)$, respectively. No
restriction on the parameter matrix $\bfq$ is needed here; it is only in
patching quantum tori together to cover a quantum affine space, as in the
following section, that we will need to avoid $-1$.

Let $k$ be a field. From (3.7) onward, we will assume that $k$ is
algebraically closed.

\definition{3.1} Let $\bfq= (q_{ij})$ be a multiplicatively
antisymmetric $n{\times}n$ matrix over $k$; that is, $q_{ii}=1$ and
$q_{ji}= q_{ij}^{-1}$, for all $i,j$. Let $A= \Oq((\kx)^n)$ be the
corresponding multiparameter quantized coordinate ring of the torus
$(\kx)^n$, that is, the $k$-algebra generated by elements
$x_1^{\pm1},\dots,x_n^{\pm1}$ subject only to the relations $x_ix_j=
q_{ij}x_jx_i$ for all $i,j$. This algebra is also known as a {\it
McConnell-Pettit algebra\/}, after \cite{\McCPet}; in the notation of
that paper, $A= P(\bfq)$.

Let us first express $A$ in terms of ordered monomials, using standard
multi-index notation. Thus $A$ has a $k$-basis of monomials,
$x^\alpha= x_1^{\alpha_1} x_2^{\alpha_2} \cdots x_n^{\alpha_n}$, for
$n$-tuples $\alpha= (\alpha_1,\dots,\alpha_n)$ from the group $\Gamma
:= \ZZ^n$.  Define $\sigma : \Gamma\times\Gamma \rightarrow \kx$ by
$$\sigma(\alpha,\beta)= \prod_{i,j=1}^n q_{ij}^{\alpha_i\beta_j}.$$
Then $\sigma$ determines the commutation rules in $A$; namely,
$x^\alpha x^\beta= \sigma(\alpha,\beta)x^\beta x^\alpha$ for
$\alpha,\beta\in \Gamma$.  Moreover, $\sigma$ is an alternating
bicharacter on $\Gamma$, because
$$\sigma (\alpha,\alpha) = 1, \qquad
\sigma(\beta,\alpha)=\sigma(\alpha,\beta)^{-1}, \qquad \text{and}
\qquad \sigma(\alpha,\beta+\beta') =
\sigma(\alpha,\beta)\sigma(\alpha,\beta'), $$
for $\alpha,\beta,\beta' \in \Gamma$.  \enddefinition

\definition{3.2} It is convenient to view $A$ as a twisted group algebra
of $\Gamma$. One way to do this is to use the function $d :
\Gamma\times\Gamma \rightarrow\kx$ such that $x^\alpha x^\beta=
d(\alpha,\beta)x^{\alpha+\beta}$ for all $\alpha,\beta\in \Gamma$. It
follows from the associativity of multiplication in $A$ that $d$ must be a
2-cocycle on
$\Gamma$. Thus by writing $A$ in terms of the $k$-basis
$\{x^\alpha \mid \alpha \in \Gamma\}$, we have expressed $A$ as a twisted
group algebra $k^d\Gamma$.

Note that $\sigma(\alpha,\beta)= d(\alpha,\beta)d(\beta,\alpha)^{-1}$ for
$\alpha,\beta\in
\Gamma$. Conversely, if $c$ is any 2-cocycle on $\Gamma$ such that
$c(\alpha,\beta) c(\beta,\alpha)^{-1}= \sigma(\alpha,\beta)$ for
$\alpha,\beta\in
\Gamma$, then
$k^c\Gamma \cong A$. Namely, $k^c\Gamma$ has a basis $\{ x_\alpha \mid
\alpha \in \Gamma \}$ such that $x_\alpha x_\beta =c(\alpha,\beta)
x_{\alpha+\beta}$ for $\alpha,\beta \in \Gamma$. If $\epsilon_1, \dots,
\epsilon_n$ denotes the standard basis for $\Gamma$, then
$k^c\Gamma$ is generated by the elements $x_{\epsilon_1}^{\pm1},
\dots, x_{\epsilon_n}^{\pm1}$, which satisfy the same commutation
rules as $x_1,\dots, x_n$, namely 
$$x_{\epsilon_i} x_{\epsilon_j}= c(\epsilon_i, \epsilon_j) x_{\epsilon_i+
\epsilon_j}= c(\epsilon_i, \epsilon_j) c(\epsilon_j, \epsilon_i)^{-1}
x_{\epsilon_j} x_{\epsilon_i}=
\sigma(\epsilon_i, \epsilon_j) x_{\epsilon_j} x_{\epsilon_i}=
q_{ij} x_{\epsilon_j} x_{\epsilon_i}$$ 
for all $i,j$. Therefore
there exists a $k$-algebra isomorphism $k^c\Gamma \rightarrow A$
such that $x_{\epsilon_i} \mapsto x_i$ for all $i$. Note that this
isomorphism sends each $x_\alpha$ to a scalar multiple of $x^\alpha$.

While the cocycle $d$ is useful for some purposes, it has some
drawbacks; for example, it is neither symmetric nor antisymmetric in
general. We therefore leave the choice of a particular cocycle until
later (see (3.5) and (3.12)).  \enddefinition

\definition{3.3} Let $c$ be any 2-cocycle on $\Gamma$ such that
$c(\alpha,\beta) c(\beta,\alpha)^{-1} =\sigma(\alpha,\beta)$ for
$\alpha,\beta \in \Gamma$, and identify $A$ with $k^c\Gamma$. This allows
us to write $A$ in terms of a $k$-basis
$\{x_\alpha \mid \alpha\in \Gamma\}$, with multiplication given by
$$x_\alpha x_\beta = c(\alpha,\beta)x_{\alpha+\beta}$$ for
$\alpha,\beta\in \Gamma$.
We shall also need the ordinary group algebra $\kgam$. Let us express
this algebra in terms of a $k$-basis $\{y_\alpha \mid \alpha\in
\Gamma\}$, with $y_\alpha y_\beta= y_{\alpha+\beta}$ for
$\alpha,\beta\in \Gamma$.  \enddefinition

\definition{3.4} Set $H= \Hom(\Gamma,\kx)$, an algebraic torus of rank
$n$ (which will henceforth be equipped with the Zariski
topology). Write the application of $H$ to $\Gamma$ in terms of a
pairing $\langle-,-\rangle : H\times \Gamma \rightarrow \kx$.  Define
actions of $H$ on both $A$ and $\kgam$ via $k$-algebra automorphisms
such that $$h{.}x_\alpha= \langle h,\alpha\rangle x_\alpha
\quad\text{and}\quad h{.}y_\alpha= \langle h,\alpha \rangle y_\alpha,$$
for $h\in H$ and $\alpha\in \Gamma$.  \enddefinition

\definition{3.5} Let $\Phi_c : A\rightarrow \kgam$ be the $k$-linear
isomorphism such that $\Phi_c(x_\alpha)= y_\alpha$ for $\alpha\in
\Gamma$. Observe that $\Phi_c$ is $H$-equivariant. The extent to which
$\Phi_c$ fails to preserve products can be expressed in terms of $c$:
since $x_\alpha x_\beta= c(\alpha,\beta) x_{\alpha+\beta}$, we have
$$\Phi_c(x_\alpha x_\beta)= c(\alpha,\beta)y_{\alpha+\beta}=
c(\alpha,\beta) \Phi_c(x_\alpha) \Phi_c(x_\beta)$$
for $\alpha,\beta\in \Gamma$.

As in \cite{\GooLet, 1.2}, set $S= \rad(\sigma)= \{\alpha \in\Gamma
\mid \sigma(\alpha,-) \equiv1\}$, a subgroup of $\Gamma$; then 
$$Z(A)=k[x_\alpha \mid \alpha\in S].$$ 
If we identify the group algebra $kS$
with the $k$-linear span of $\{y_\alpha \mid \alpha\in S\}$ inside
$k\Gamma$, then
$\Phi_c(Z(A))= kS$. In fact: \enddefinition

\proclaim{Lemma} Assume that $c\equiv 1$ on $S\times \Gamma \subseteq
\Gamma \times \Gamma$. Then the map $\Phi_c$ restricts to a
$k$-algebra isomorphism of $Z(A)$ onto $kS$. With respect to this
isomorphism, $\Phi_c$ itself gives a semilinear isomorphism of $A$ as
$Z(A)$-module onto $\kgam$ as $kS$-module: $\Phi_c(za)= \Phi_c(z)
\Phi_c(a)$ for all $z\in Z(A)$ and $a\in A$.\endproclaim

\demo{Proof} By assumption,
$\Phi_c(x_\alpha x_\beta)= \Phi_c(x_\alpha) \Phi_c(x_\beta)$ for 
$\alpha \in S$ and $\beta \in \Gamma$. Consequently, $\Phi_c(za)= \Phi_c(z)
\Phi_c(a)$ for all
$z\in Z(A)$ and $a\in A$, and the lemma follows. \qed\enddemo

We shall need the following fact to show that $S$ equals the
intersection of the kernels of some homomorphisms from $H$ (see
(3.7)).

\proclaim{3.6 Lemma} The order of the torsion subgroup of
$\Gamma/S$ is not divisible by $\chr k$.\endproclaim

\demo{Proof} The torsion subgroup
of $\Gamma/S$ is a finite abelian group, hence a finite direct
sum of finite cyclic groups, and so the order of this subgroup
equals the product of the orders of certain of its elements. Thus it
suffices to show that the order of each torsion element of
$\Gamma/S$ is not divisible by $\chr k$.

Let $\alpha+S$ be a torsion element of $\Gamma/S$, and choose a
basis $\epsilon_1 ,\dots, \epsilon_n$ for $\Gamma$. For each positive
integer $m$, we have
$$m\alpha\in S \quad\Longleftrightarrow\quad \sigma(m\alpha,
\gamma_i)=1 \text{\ for all\ } i \quad\Longleftrightarrow\quad
\sigma(\alpha,
\epsilon_i)^m=1 \text{\ for all\ } i.$$
Hence, the order of $\alpha+S$ equals the order of the element
$(\sigma(\alpha, \epsilon_1), \dots, \sigma(\alpha, \epsilon_n))$ in
$(\kx)^n$, and the latter order is clearly not divisible by $\chr
k$. \qed\enddemo

\definition{3.7} From now until the end of the section, assume that
$k$ is algebraically closed.

We use
$^\perp$ to denote orthogonals of subsets of
$H$ and $\Gamma$, relative to the pairing $\langle-,-\rangle$. In
particular,
$$\sperp= \{h\in H\mid \langle h,\alpha\rangle=1 \ \text{for} \
\alpha\in S\}= \{h\in H\mid \ker h\supseteq S\},$$
a closed subgroup of $H$, and
$$S^{\perp\perp}= \{\alpha \in\Gamma \mid \langle h,\alpha\rangle=1
\ \text{for} \ h\in \sperp\}= \bigcap_{h\in \sperp} \ker h,$$
a subgroup of $\Gamma$.
\enddefinition

\proclaim{Lemma} $S^{\perp\perp}= S$.\endproclaim

\demo{Proof} Consider $\alpha\in \Gamma\setminus S$. The order of
the coset $\alpha +S$ in $\Gamma/S$ is either infinite, or -- if
finite -- not divisible by $\chr k$, in view of (3.6). In either
case, the subgroup $\ZZ(\alpha +S)\subseteq \Gamma/S$ supports a
nontrivial homomorphism to
$\kx$. Since
$\kx$ is divisible, any such homomorphism extends to a homomorphism
$\Gamma/S
\rightarrow \kx$. Thus, there exists a homomorphism $h : \Gamma
\rightarrow \kx$ (i.e., $h\in H$) such that $S\subseteq \ker h$ and
$h(\alpha) \ne1$, which verifies that $\alpha \notin S^{\perp\perp}$.
\qed\enddemo

\proclaim{3.8 Proposition} {\rm (a)} $A^{\sperp}= Z(A)$ and
$(\kgam)^{\sperp}= kS$. 

{\rm (b)} Every $\sperp$-stable ideal of
$\kgam$ is generated by its contraction to $kS$.

{\rm (c)} Every ideal of $A$ is generated by its contraction to $Z(A)$.

{\rm (d)} Assume that $c\equiv 1$ on $S\times \Gamma$. The map $I\mapsto
\Phi_c(I)$, for ideals
$I$ of
$A$, provides a product preserving, $H$-equivariant lattice isomorphism
from the lattice of ideals of $A$ onto the lattice of $\sperp$-ideals of
$\kgam$.  \endproclaim

\demo{Proof} (a) It follows from (3.7) that $A^{\sperp}$ is spanned
by the $x_\alpha$ for $\alpha\in S$, and that $(\kgam)^{\sperp}$
is spanned by the $y_\alpha$ for $\alpha\in S$. The second part of (a)
is thus clear, and the first follows from (3.5).

(b) Consider an $\sperp$-stable ideal $I$ in $\kgam$. Since $\sperp$
acts semisimply on $\kgam$, this ideal must be spanned by
$\sperp$-eigenvectors. Let $a\in I$ be an $\sperp$-eigenvector, and
choose an element $\gamma$ from the support of $a$. Then
$x_\gamma^{-1}a$ is also an $\sperp$-eigenvector and contains $1$ in
its support. Hence $\sperp$ fixes $x_\gamma^{-1}a$. Thus
$x_\gamma^{-1}a \in (\kgam)^{\sperp}= kS$, and so $a\in (\kgam)(I\cap
kS)$. Therefore $I= (\kgam)(I\cap kS)$.

(c) See, for example, \cite{\GooLet, 1.4}.

(d) To start, note that $U \mapsto \Phi_c(U)$ produces an
$H$-equivariant, bijective lattice isomorphism from the lattice of
$k$-subspaces $U$ of $A$ onto the lattice of $k$-subspaces of $\kgam$.

Next, let $I$ be an ideal of $A$. It follows from (c) that $I = I'A$
for some ideal $I'$ of $Z(A)$. By (3.5), $\Phi_c(I')$ is an ideal of
$kS$, and $\Phi_c(I) = \Phi_c(I')\Phi_c(A) = \Phi_c(I')\kgam$ is an
ideal of $\kgam$. Furthermore, since $\Phi_c(I')$ is an ideal of $kS$,
and since $kS$ is fixed pointwise by $\sperp$, we see that
$\Phi_c(I')\kgam$ is stable under $\sperp$. Similarly, it follows from
(b) and (3.5) that $\Phi_c^{-1}$ maps each $\sperp$-ideal of $\kgam$
to an ideal of $A$.  Therefore the set maps $\Phi_c$ and $\Phi_c^{-1}$
provide mutually inverse bijections between the lattices of ideals of
$A$ and $\sperp$-ideals of $k\Gamma$. It follows easily from the
preceding paragraph that these bijections are $H$-equivariant lattice
isomorphisms.

Now let $I$ and $J$ be ideals of $A$. It only remains to show that
$\Phi_c(IJ) = \Phi_c(I)\Phi_c(J)$. As above, write $I = I'A$ and $J =
J'A$, for ideals $I'$ and $J'$ of $Z(A)$. Then $IJ=I'J'A$, and so by
(3.5), $\Phi_c(IJ) = \Phi_c(I')\Phi_c(J')\kgam = \Phi_c(I)\Phi_c(J)$.
\qed\enddemo

\definition{3.9} Let $\spspec\kgam$ and $\spmax\kgam$ denote the sets
of $\sperp$-prime and maximal proper $\sperp$-ideals of $\kgam$,
respectively. As noted in (1.11), $\spmax\kgam \subseteq
\spspec\kgam$. Equip these sets with their respective Zariski
topologies, again following (1.11).
\enddefinition

\proclaim{Proposition} Assume that $c\equiv 1$ on $S\times \Gamma$.  The
set maps
$\Phi_c$ and
$\Phi_c^{-1}$ provide pairs of mutually inverse, $H$-equi\-var\-iant,
homeomorphisms
$$\spec A\longleftrightarrow \spspec\kgam \qquad\text{and}\qquad \max
A\longleftrightarrow \spmax\kgam.$$
\endproclaim

\demo{Proof} This follows from (3.8d). \qed\enddemo

\proclaim{3.10 Proposition} The assignment $P \mapsto (P : \sperp)$,
for prime ideals $P$ of $\kgam$, produces $H$-equivariant topological
quotient maps
$$\spec \kgam \rightarrow \spspec \kgam \qquad \text{and}
\qquad \max \kgam \rightarrow \spmax \kgam .$$ \endproclaim

\demo{Proof} By (1.13a), $\spmax \kgam= \{ (M:\sperp) \mid M\in\max
\kgam \}$, and so the given assignment does map $\max\kgam$ onto $\spmax
\kgam$. The given maps are
$H$-equivariant because
$H$ is abelian. The remainder of the proposition follows from (1.11).
\qed\enddemo

The two propositions above provide the information necessary to exhibit $\spec
A$ and $\prim A$ as topological quotients of their classical analogs, as in
the following theorem. A cocycle $c$ with the required condition is
constructed in (3.12).

\proclaim{3.11 Theorem} Assume that $c\equiv 1$ on $S\times \Gamma$. Then
the rule $P
\mapsto
\Phi_c^{-1}(P:\sperp)$, for prime ideals $P$ of $\kgam$, defines
$H$-equivariant topological quotient maps
 $$\spec\kgam
\rightarrow \spec A \qquad\text{and}\qquad \max\kgam
\rightarrow \max A.$$
The fibers of the second map are exactly the $\sperp$-orbits in
$\max\kgam$. \endproclaim

\demo{Proof} The first part follows from (3.9) and (3.10).
The final part is a consequence of (1.13b). \qed\enddemo

\definition{3.12} In order to apply (3.11), we must be able to choose the
cocycle $c$ to be trivial on $S\times \Gamma$. This can be done as follows:
\enddefinition

\proclaim{Lemma} There exists a $2$-cocycle $c$ on $\Gamma$ such that
$c\equiv 1$ on $S\times \Gamma$ and $c(\alpha,\beta) c(\beta,\alpha)^{-1}
= \sigma(\alpha,\beta)$ for $\alpha,\beta \in \Gamma$. \endproclaim

\demo{Proof} As in \cite{\GooLet, 1.6}, $\Gamma$ has a basis
$\gamma_1, \dots,  \gamma_n$ such that $S$ is generated by
$m_1\gamma_1 ,\dots, m_t \gamma_t$, for some $t\le n$ and some positive
integers $m_i$. Set $p_{ij}= \sigma( \gamma_i, \gamma_j)$ for all
$i,j$, and note that the matrix $(p_{ij})$ is multiplicatively
antisymmetric. Given $\alpha, \beta \in \Gamma$, write $\alpha= \sum_i
a_i\gamma_i$ and $\beta= \sum_j b_j\gamma_j$ for some $a_i,b_j \in
\ZZ$; then $$\sigma(\alpha, \beta)= \prod_{i,j} p_{ij}^{a_ib_j}.$$ Now
define $$c(\alpha, \beta)= \prod_{i<j} p_{ij}^{a_ib_j}$$ for $\alpha=
\sum_i a_i\gamma_i$ and $\beta= \sum_j b_j\gamma_j$ in $\Gamma$. Then
$c$ is a bicharacter on $\Gamma$, and hence a 2-cocycle. It is easily
checked that $c(\alpha, \beta) c(\beta, \alpha)^{-1} = \sigma(\alpha,
\beta)$ for $\alpha, \beta \in \Gamma$.  Since $m_i\gamma_i \in S$, we
obtain $$1= \sigma(m_i\gamma_i, \gamma_j) =p_{ij}^{m_i}$$ for all
$i,j$. Hence, if $\alpha= \sum_i a_i\gamma_i \in S$, then
$p_{ij}^{a_i} =1$ for all $i,j$, and therefore $c(\alpha, \beta) =1$
for all $\beta \in \Gamma$.  \qed\enddemo

\definition{3.13 Remark} It would result in a stronger theorem if we
could conclude in (3.11) that the given maps were closed. However,
this conclusion is not possible in general. For instance, take $n=3$,
choose $p\in \kx$ not a root of unity, and take $$\bfq= \left(
\smallmatrix 1&1&1\\ 1&1&p^2\\ 1&p^{-2}&1 \endsmallmatrix \right).$$
Then $\sigma$ is given by the rule $\sigma(s,t)=
p^{2(s_2t_3-s_3t_2)}$, and we may choose $c$ to be given by $$c(s,t)=
p^{s_2t_3-s_3t_2}.$$ We compute that $S= \ZZ(1,0,0)$ and $\sperp= \{
(h_1,h_2,h_3) \in H\mid h_1 =1\}$.

Now write $k\Gamma= k[y_1^{\pm1}, y_2^{\pm1}, y_3^{\pm1}]$, where
$y_1,y_2,y_3$ correspond to the standard basis for $\Gamma$, and let
$$V= \{ \langle y_1-\alpha_1,\, y_2-\alpha_2,\, y_3-\alpha_3 \rangle \mid
\alpha_i \in \kx \text{\ and\ } (\alpha_1-1)\alpha_2=1 \},$$ 
a closed
subset of $\max\kgam$. The map $M\mapsto (M:\sperp)$ sends $V$ to the set
$$\{ \langle y_1-\alpha_1\rangle \mid \alpha_1\in\kx \text{\ and\ }
\alpha_1\ne1 \},$$ 
a non-closed subset of $\spmax\kgam$, and thus $M\mapsto
\Phi_c^{-1}(M:\sperp)$ sends $V$ to a non-closed subset of $\max A$.
\enddefinition

\head 4. Quantum affine spaces\endhead

We now present and prove the main results of this paper. In essence, we
partition the prime and primitive spectra of a quantum affine space into
corresponding spectra of quantum tori, apply the results of the previous
section to these spectra, and then patch everything together. It is the
patching process that requires a more careful choice of cocycle than
previously.

We continue
to assume that the base field $k$ is algebraically closed.

\definition{4.1} Again let $\bfq= (q_{ij})$ be a multiplicatively
antisymmetric $n\times n$ matrix over $k$. Now, however, let $A=
\Oq(k^n)$ be the corresponding multiparameter quantized coordinate
ring of affine $n$-space. In other words, $A$ is the $k$-algebra
generated by elements $x_1,\dots,x_n$, subject only to the relations
$x_ix_j= q_{ij}x_jx_i$. We will treat $A$ as a twisted semigroup
algebra, as follows.

Set $\Gamma =\ZZ^n$ and $\Gamma^+= (\ZZ^+)^n$, and define an
alternating bicharacter $\sigma$ on $\Gamma$ by the rule
$$\sigma(\alpha,\beta)= \prod_{i,j=1}^n q_{ij}^{\alpha_i \beta_j}$$
as in (3.1). If $c$ is any 2-cocycle on $\Gamma$ such that
$c(\alpha,\beta) c(\beta,\alpha)^{-1}= \sigma(\alpha,\beta)$ for
$\alpha,\beta \in \Gamma$, we can write
$A$ as the twisted semigroup algebra $k^c\Gamma^+$, viewed in the
obvious way as a subalgebra of the twisted group algebra
$k^c\Gamma$. 

We will need to apply the results of Section 3 to the localizations
$$A_w= \bigl( A/\langle x_i\mid i\in w\rangle \bigr) [x_j^{-1}\mid j\notin
w]$$
as $w$ ranges over the subsets of $\{1,\dots,n\}$. Thus, we require a
cocycle $c$ whose restrictions to the subgroups of $\Gamma$ generated by
subsets of the standard basis all satisfy the hypothesis of (3.11). The
cocycle obtained in (3.12) is not sufficient for this purpose. With a
minor technical assumption, we can construct a suitable cocycle as follows.
\enddefinition

\definition{4.2} Let $\langle q_{ij}\rangle$ denote the subgroup of $\kx$
generated by the $q_{ij}$. Assume from now on that either
$-1\notin \langle q_{ij}\rangle$ or  $\chr k=2$.
\enddefinition

\proclaim{Lemma} There exists an alternating bicharacter $c$ on $\Gamma$
such that $c^2=\sigma$ and such that $c(\alpha,\beta)=1$ whenever
$\sigma(\alpha,\beta)=1$. \endproclaim

\demo{Proof} Since $k$ is algebraically closed, the abelian group
$\kx$ is divisible, and so it contains a divisible hull $D$ for the
subgroup $\langle q_{ij}\rangle$. (The divisible hull is the injective
hull within the category of abelian groups.) In view of our
hypotheses, $\langle q_{ij}\rangle$ contains no elements of order 2,
whence the same holds for $D$. Now each $q_{ij}$ has a unique square
root in $D$, say $p_{ij}$, and it follows from the uniqueness of the
choices that $(p_{ij})$ is a multiplicatively antisymmetric matrix
over $k$. Thus, the rule $$c(\alpha,\beta)= \prod_{i,j=1}^n
p_{ij}^{\alpha_i \beta_j}$$ defines an alternating bicharacter $c$ on
$\Gamma$ such that $c^2=\sigma$.  For any $\alpha,\beta \in \Gamma$
such that $\sigma(\alpha,\beta)=1$, we have $c(\alpha,\beta)^2=1$ and
therefore $c(\alpha,\beta)=1$, because $D$ has no elements of order
2. \qed\enddemo

\definition{4.3} Fix an alternating bicharacter $c$ as in (4.2). Then $c$
is a 2-cocycle on $\Gamma$ and 
$$c(\alpha,\beta) c(\beta,\alpha)^{-1}=
c(\alpha,\beta)^2= 
\sigma(\alpha,\beta)$$
for
$\alpha,\beta \in \Gamma$. As in (3.3), we write $k^c\Gamma$ in terms of a
basis
$\{x_\alpha\mid \alpha\in \Gamma\}$ such that $x_\alpha x_\beta=
c(\alpha,\beta) x_{\alpha+\beta}$ for $\alpha,\beta \in \Gamma$.
Thus $A= k^c\Gamma^+$ equals the subspace of $k^c\Gamma$ spanned by
the
$x_\alpha$ for $\alpha \in \Gamma^+$.

Similarly, we write the group algebra $\kgam$ in terms of a basis
$\{y_\alpha \mid \alpha\in \Gamma\}$, with $y_\alpha y_\beta=
y_{\alpha+\beta}$ for $\alpha,\beta\in \Gamma$, and we identify the
semigroup algebra $R:= \kgam^+$ with the subspace of $\kgam$ spanned by
the $y_\alpha$ for $\alpha \in \Gamma^+$. We also view $R$ as a
polynomial ring $k[y_1,\dots, y_n]$, as in Section 2, where $y_i=
y_{\epsilon_i}$ and $\epsilon_i$ is the $i$-th standard basis element
of $\Gamma$.  \enddefinition

\definition{4.4} Following (3.4), set $H= \Hom(\Gamma,\kx)$ and write
the application of $H$ to $\Gamma$ in terms of a pairing
$\langle-,-\rangle : H\times \Gamma \rightarrow \kx$. We again have
actions of $H$ on $k^c\Gamma$ and $\kgam$ via $k$-algebra
automorphisms such that $h{.}x_\alpha= \langle h,\alpha\rangle
x_\alpha$ and $h{.}y_\alpha= \langle h,\alpha \rangle y_\alpha$, for
$h\in H$ and $\alpha\in \Gamma$.  The subalgebras $A$ and $R$ are stable
under these actions.  \enddefinition

\definition{4.5} Let
$W$ denote the set of subsets of $\{1,\ldots,n\}$. For $w\in W$, let
$\spec_w A$ be the set, equipped with the relative Zariski topology, of
those prime ideals $P$ in $A$ such that $$P\cap \{x_1,\dots,x_n\} = \{x_i
\mid i\in w\}.$$
 Then $\spec A$ is the disjoint union of the sets $\spec_w A$.
Likewise, $\prim A$ is the disjoint union of the sets $\prim_w A=
\prim A\cap \spec_w A$, each of which is also endowed with the
relative topology. Note that the sets $\spec_w A$ and $\prim_w A$ are
invariant under the action of $H$. Also, it was proved in
\cite{\GooLet, 2.3}, for $w \in W$, that the set $\prim _w A$
coincides with the set of maximal members of $\spec _w A$. Define
$\spec_w R$ and $\max_w R$ similarly, as in (2.1) and (2.5). Recall
that $\max_w R$ equals the set of maximal elements of $\spec_w R$.

It is convenient to label the localizations 
$$A_w= \bigl( A/\langle
x_i\mid i\in w\rangle \bigr) [x_j^{-1}\mid j\notin w]
\qquad\text{and}\qquad R_w= \bigl( R/\langle y_i\mid i\in w\rangle
\bigr) [y_j^{-1}\mid j\notin w],$$
and to identify $A_w$ (respectively, $R_w$) with the $k$-subalgebra of
$k^c\Gamma$ (respectively, $k\Gamma$) spanned by the $x_\alpha$
(respectively, $y_\alpha$) for those $\alpha\in \Gamma$ such that
$\alpha_i=0$ for $i\in w$.  In view of the preceding paragraph and
\cite{\GooWar, 9.20}, for example, we see that localization induces
$H$-equivariant homeomorphisms $$\alignat2 \spec_w A &\longrightarrow
\spec A_w &\qquad\qquad \spec_w R &\longrightarrow \spec R_w\\ \prim_w
A &\longrightarrow \max A_w &\qquad\qquad \max_w R &\longrightarrow
\max R_w. \endalignat$$ \enddefinition

\definition{4.6} As in (3.5), we have an $H$-equivariant $k$-linear
isomorphism $\Phi =\Phi_c : k^c\Gamma \rightarrow \kgam$ such that
$\Phi(x_\alpha) =y_\alpha$ for $\alpha\in \Gamma$. This map restricts
to an $H$-equivariant $k$-linear isomorphism from $A$ onto $R$. Note that
$\Phi$ also restricts to $H$-equivariant $k$-linear
isomorphisms $\Phi_w : A_w\rightarrow R_w$ for $w \in W$.
\enddefinition

\definition{4.7} Retain $\epsilon_1, \dots, \epsilon_n$ as the
standard basis for $\Gamma$, and fix $w\in W$. Let $\Gamma_w$ be the
subgroup of $\Gamma$ generated by $\{ \epsilon_i \mid i\notin w\}$,
and let $\sigma_w$ and $c_w$ denote the restrictions of $\sigma$ and
$c$ to $\Gamma_w$. With this notation, the identifications made in
(4.5) are $A_w= k^{c_w}\Gamma_w$ and $R_w= k\Gamma_w$.  \enddefinition

\definition{4.8} Set $S_w= \rad(\sigma_w)= \{ \alpha\in \Gamma_w \mid
\sigma(\alpha,-) \equiv1 \text{\ on\ } \Gamma_w \}$, for $w\in W$. As
in (3.7), set
$$S_w^\perp = \left\{ h \in H \mid \text{$\langle h , \alpha \rangle =
1$ for $\alpha \in S_w$}\right\} \quad \text{and} \quad
\Gamma_w^\perp = \left\{ h \in H \mid \text{$\langle h , \alpha
\rangle = 1$ for $\alpha \in \Gamma_w$}\right\}. $$
Because $\sigma \equiv 1$ on $S_w\times \Gamma_w$, it follows from
(4.2) that $c\equiv 1$ on $S_w\times \Gamma_w$.  \enddefinition

\proclaim{Lemma} Let $v\subseteq w$ in $W$.

{\rm (a)} $\swperp \subseteq S_v^\perp \Gamma_w^\perp$.

{\rm (b)} $(P:S_v^\perp)\subseteq (P:\swperp)$ for all $P\in
\spec_w R$. \endproclaim

\demo{Proof} (a) Observe that $\Gamma_w\subseteq \Gamma_v$ and that
$S_v\cap \Gamma_w\subseteq S_w$. Any $h\in \swperp$ is a
homomorphism $\Gamma \rightarrow \kx$ such that $\ker h\supseteq
S_w\supseteq S_v\cap \Gamma_w$. Hence, $h|_{\Gamma_w}$ extends to a
homomorphism $h_1 : \Gamma_w +S_v\rightarrow \kx$ such that
$S_v\subseteq \ker h_1$. Since $\kx$ is divisible, $h_1$ extends to a
homomorphism $h_2 : \Gamma\rightarrow \kx$. Now $h_2$ is a
homomorphism in $H$ such that $h_2=h$ on $\Gamma_w$ and $h_2=1$ on
$S_v$. Thus $h_2\in S_v^\perp$ and $h_2^{-1}h\in \Gamma_w^\perp$, whence
$h\in S_v^\perp \Gamma_w^\perp$ as desired.

(b) Observe that $\Gamma_w^\perp$ acts trivially on $R/\langle y_i \mid
i\in w \rangle$. Hence, $P$ is stable under $\Gamma_w^\perp$, and so
$(P:S_v^\perp)= (P:S_v^\perp
\Gamma_w^\perp)$. In view of part (a), $(P:S_v^\perp \Gamma_w^\perp)
\subseteq (P:\swperp)$, and we are done.  \qed\enddemo

\definition{4.9} Set $\G = \{ \sperp _w \mid w \in W \}$, and note by
(4.8) that $\G$ satisfies the hypotheses in (2.2), as applied to
$R$. Define $\ggspec R$ as in (2.2), and define $(P : \G)$ to be $(P :
\sperp_w)$ for $P \in \spec _w R$, as in (2.3). Apply the Zariski
topology to $\ggspec R$, following (2.4). We may also use (2.4) to
conclude that the assignment $P \mapsto (P : \G)$ is a topological
quotient map from $\spec R$ onto $\ggspec R$.

Next, define $\ggmax R$ as in (2.5). It follows from (2.6), for each
$w \in W$, that $\ggmax _wR$ is the set of maximal members of $\ggspec
_wR$. From (2.5) it follows that the assignment $M \mapsto (M : \G)$,
for $M \in \max R$, produces an $H$-equivariant topological quotient
map from $\max R$ onto $\ggmax R$. Moreover, by (2.7), the fibers in
$\max R$ over points in $\ggmax _w R$ are precisely the $\sperp
_w$-orbits in $\max _w R$.  \enddefinition

It remains to relate $\spec A$ to $\ggspec R$ and $\prim A$ to $\ggmax R$.

\proclaim{4.10 Lemma} The set function $I \mapsto \Phi(I)$, for ideals $I$
of $R$, produces an $H$-equivariant homeomorphism from $\spec A$ onto
$\ggspec R$ that restricts to a homeomorphism from
$\prim A$ onto $\ggmax R$. \endproclaim

\demo{Proof} Given $w\in W$, there exists an $H$-equivariant
commutative diagram
$$\CD A @>\Phi>> R\\ @V\lambda_w VV @VV\mu_w V\\ A_w @>\Phi_w>> R_w
\endCD$$
where $\lambda_w$ and $\mu_w$ are the respective localization
maps. Let $Q\in \ggspec_w R$, and write $Q=(P:\swperp)$ for some $P\in
\spec_w R$.  Then $P= \mu_w^{-1}(P')$ for some $P'\in \spec R_w$; the
ideal $Q'= (P':\swperp)$ lies in $\swpspec R_w$; and $Q=
\mu_w^{-1}(Q')$. By (3.9), $\Phi_w^{-1}(Q') \in \spec A_w$, and so
$\lambda_w^{-1}\Phi_w^{-1}(Q') \in \spec_w A$. Since
$\lambda_w^{-1}\Phi_w^{-1}(Q') =\Phi^{-1}\mu_w^{-1}(Q')=
\Phi^{-1}(Q)$, we thus see that the set $\Phi^{-1}(Q)$ is a prime
ideal of $A$, lying in $\spec_w A$. Similarly, $\Phi^{-1}(M)\in
\prim_w A$ for all $M\in \ggmax_w R$.

Since $\Phi$ and $\Phi_w$ are bijections, analogous arguments show that
$\Phi(Q)\in \ggspec_w R$, for all $Q\in \spec_w A$, and $\Phi(P)\in \ggmax_w
R$, for all $P\in \prim_w A$.

The results above, taken over all $w\in W$, show that the set maps
$I\mapsto \Phi(I)$ and $J\mapsto \Phi^{-1}(J)$ yield pairs of mutually
inverse bijections $$\spec A \longleftrightarrow \ggspec R
\qquad\text{and}\qquad \prim A \longleftrightarrow \ggmax R.$$ These
maps are $H$-equivariant because $\Phi$ and $\Phi^{-1}$ are
$H$-equivariant. Finally, note that the topologies on all of the
spaces considered in this section can be defined using subsets of $A$
and $R$ instead of ideals -- for instance, the closed subsets of
$\spec A$ are precisely those of the form $\{Q\in \spec A \mid
Q\supseteq X\}$, for arbitrary subsets $X \subseteq A$. It therefore
follows that the above maps are homeomorphisms, as
desired. \qed\enddemo

Combining (4.9) and (4.10), we obtain the main theorem of the paper,
which we state as follows, using the notation developed in
(4.2)--(4.8).

\proclaim{4.11 Theorem} Suppose that $k$ is an algebraically closed
field, that $\bfq= (q_{ij})$ is a multiplicatively antisymmetric
$n\times n$ matrix over $k$, that $A= \Oq(k^n)$, and that $R=
\O(k^n)$. Further assume that $-1 \notin \langle q_{ij}\rangle$ or
that $\chr k=2$. Let $H=(\kx)^n$ act by $k$-algebra automorphisms on $A$
 and on $R$ in the standard manner. Then there exist
$H$-equivariant topological quotient maps
$$\spec R\rightarrow \spec A \qquad\text{and}\qquad \max R \rightarrow
\prim A$$
given by $P\mapsto \Phi^{-1}(P:\swperp)$, for $P\in \spec_w R$. The
fibers of the second map, over points in $\prim_w A$, consist precisely
of the $\swperp$-orbits within $\max_w R$. \qed\endproclaim

\head 5. Some examples\endhead

We illustrate (4.11) by calculating the explicit form of the maps from $\max
R$ to $\prim A$ in the cases of the standard single parameter quantum affine
spaces, and we comment on some more general cases related to bilinear forms
and Poisson brackets.
Let $k$, $\bfq$, $A$, and $R$ be as in (4.11), and let $$\Psi : \max R
\longrightarrow \prim A$$ denote the topological quotient map given by
the theorem therein. Composing $\Psi$ with the natural isomorphism of
affine $n$-space onto $\max R$, we obtain a similar topological
quotient map 
$$\psi : k^n \longrightarrow \prim A.$$ 
For simplicity of notation, we describe $\psi$ rather than $\Psi$.

\definition{5.1} First, let $A$ be the standard one-parameter quantization
$\O_q(k^n)$, for some $q\in\kx$. Then $A$ is generated by $x_1,\dots,x_n$
such that $x_ix_j= qx_jx_i$ for all $i<j$. The alternating bicharacter
$\sigma$ as in (4.1) can be expressed by
$$\sigma(\alpha,\beta) = q^{b(\alpha,\beta)},$$
where $b : \Gamma\times \Gamma \rightarrow \ZZ$ is the alternating
bilinear form 
$$b(\alpha, \beta) = \sum_{i<j} \alpha_i\beta_j -\sum_{i>j} \alpha_i
\beta_j.$$

To complete the hypotheses of (4.11), we must assume that $q$ is either
not a root of unity or an odd root of unity. In either case, $q$ has a
square root $p\in\kx$ such that $-1\notin \langle p\rangle$ (if $\chr k\ne
2$), and the rule
$$c(\alpha,\beta) = p^{b(\alpha,\beta)}$$
defines an alternating bicharacter $c$ on $\Gamma$ satisfying the
conclusions of (4.2). Hence, we can identify $A$ with $k^c\Gamma^+$ for
this $c$.
\enddefinition

\definition{5.2} Continue with (5.1), and suppose that $q$ is not a root
of unity. We calculate $\psi(\bflam)$ for points $\bflam= (\lambda_1,
\dots, \lambda_n)$ in $k^n$. This depends on the position of $\bflam$, with
respect to the stratification of $k^n$ corresponding to the subsets
$\max_w R$. Set
$$(k^n)_w = \{ \bflam \in k^n \mid \lambda_i = 0 \text{\ for\ } i\in
w \text{\ and\ } \lambda_j \ne 0 \text{\ for\ } j\notin w \}$$
for $w\in W$.

If $n-|w|$ is even, we calculate that $S_w=0$. In this case,
$$\psi(\bflam)= \langle x_i \mid i\in w\rangle$$
for all $\bflam\in (k^n)_w$.

Now suppose that $n-|w|$ is odd. List the elements of the complement of
$w$ in ascending order, say
$$\{1,\dots,n\} \setminus w= \{ w_1 < w_2 <\cdots< w_{2m+1} \},$$
and then set
$$w_+= \{ w_1, w_3, \dots, w_{2m+1}\} \qquad\text{and}\qquad w_-= \{w_2,
w_4, \dots, w_{2m} \}.$$
Finally, write $\gamma_{\pm}(w)= \sum_{i\in w_{\pm}} \epsilon_i \in
\Gamma$, where $\epsilon_1,\dots, \epsilon_n$ is, as usual, the standard
basis for
$\Gamma$. We calculate that $S_w= \ZZ \bigl( \gamma_+(w)- \gamma_-(w) \bigr)$,
and that
$$\psi(\bflam)= \langle x_i\mid i\in w\rangle + \bigl\langle \bigl(
\prod_{j\in w_-} \lambda_j \bigr) x_{\gamma_+(w)} -\bigl( \prod_{j\in w_+}
\lambda_j \bigr) x_{\gamma_-(w)} \bigr\rangle$$
for $\bflam\in (k^n)_w$.

If one wishes to write the last formula in terms of ordinary monomials in
the $x_i$, additional scalar factors are introduced. For instance, suppose
that $n=3$ and $w=\varnothing$, so that $w_+=\{1,3\}$ and $w_-= \{2\}$.
Since $x_1x_3= c(\epsilon_1,\epsilon_3) x_{(1,0,1)}= px_{\gamma_+(w)}$, we
have
$$\psi(\lambda_1,\lambda_2,\lambda_3)= \langle p^{-1}\lambda_2
x_1x_3 - \lambda_1\lambda_3 x_2 \rangle$$
for $(\lambda_1,\lambda_2,\lambda_3) \in (k^3)_w= (\kx)^3$.
\enddefinition

\definition{5.3} Continue with (5.1), but suppose now that $q$ is a
primitive $t$-th root of unity, for some odd $t>1$. In this case, the
correct choice for $p$ is $q^{(t+1)/2}$. The form of $\psi(\bflam)$ for
$\bflam\in (k^n)_w$ again depends on the parity of $n-|w|$.

If $n-|w|$ is even, we calculate that $S_w=t\Gamma_w$, and that
$$\psi(\bflam)= \langle x_i \mid i\in w\rangle+ \langle x_j^t-
\lambda_j^t \mid j\notin w \rangle$$ 
for $\bflam\in (k^n)_w$.

If $n-|w|$ is odd, we calculate that $S_w=t\Gamma_w+ \ZZ \bigl( \gamma_+(w)-
\gamma_-(w) \bigr)$, and that
$$\psi(\bflam)= \langle x_i\mid i\in w\rangle + \langle x_j^t-
\lambda_j^t \mid j\notin w \rangle + \bigl\langle \bigl(
\prod_{j\in w_-} \lambda_j \bigr) x_{\gamma_+(w)} -\bigl( \prod_{j\in w_+}
\lambda_j \bigr) x_{\gamma_-(w)} \bigr\rangle$$
for $\bflam\in (k^n)_w$.
\enddefinition

\definition{5.4} Let $b : \Gamma\times \Gamma \rightarrow \ZZ$ be an
arbitrary alternating bilinear form, let $q\in\kx$ be a non-root of
unity, and set $q_{ij}= q^{b(\epsilon_i, \epsilon_j)}$ for all $i,j$. The
alternating bicharacter $\sigma$ associated with $A$ as in (4.1)
then takes the form $\sigma(\alpha, \beta)= q^{b(\alpha, \beta)}$. To
obtain a bicharacter $c$ as in (4.2), take a square root $p$ for $q$, and
define $c(\alpha, \beta)= p^{b(\alpha, \beta)}$. Observe that since $q$ is
not a root of unity, $S_w= \rad(\sigma|_{\Gamma_w})=
\rad(b|_{\Gamma_w})$ for $w\in W$.

Now suppose that $\chr k=0$. There is a Poisson bracket on $R$
such that $\{y_\alpha, y_\beta \}= b(\alpha, \beta) y_\alpha y_\beta$ for
$\alpha, \beta \in \Gamma$. The analysis in \cite{\Van, 1.4} (with $\CC$
replaced by $k$) shows that the $\swperp$-orbits in $(k^n)_w$ coincide
with the minimal Poisson subvarieties relative to the Poisson structure
just defined. Therefore the fibers of $\psi$ in this case coincide with
the minimal Poisson subvarieties of $k^n$.

When $k=\CC$, these minimal Poisson subvarieties are not necessarily the
same as the symplectic leaves in $k^n$. See \cite{\Van, 1.4} for necessary
and sufficient conditions. 
\enddefinition

\definition{5.5} The relation $\sigma= q^b$ in (5.4) is stronger than
necessary for the analysis given -- it would suffice to have an
alternating bilinear form $b : \Gamma\times \Gamma \rightarrow k$ such
that $\rad(\sigma|_{\Gamma_w})= \rad(b|_{\Gamma_w})$ for all $w\in
W$. To find further examples where such $b$ exist, and some situations
in which the fibers of a map analogous to $\psi$ coincide with
symplectic leaves, see \cite{\Van, Section 3}.  \enddefinition

\definition{5.6} To test whether the hypothesis $-1\notin \langle q_{ij}
\rangle$ in (4.11) is necessary, the natural examples to investigate are
the algebras $\O_{-1}(k^n)$ in characteristic different from
2. Our method works precisely when
$\O_{-1}(k^n)$ can be expressed in the form $k^c\Gamma^+$ for some
2-cocycle
$c$ such that
$c\equiv 1$ on
$S_w\times \Gamma_w$ for all $w\in W$. It is easy to find such a $c$ in
case $n=2$, and one also exists in case $n=3$. However, it can be shown
that no such $c$ exists when $n\ge 4$. Thus, our method cannot handle
$\Oq(k^n)$ for arbitrary $\bfq$ and $n$. We leave it as an open question
whether $\prim \Oq(k^n)$ is always an $H$-equivariant topological quotient
of $\max \O(k^n)$.
\enddefinition

\head 6. Twists of graded commutative algebras  \endhead

To conclude the paper, we apply the main theorem (4.11) to a class of
twists of finitely generated commutative graded algebras. Assume
throughout that $k$ is a field, that $R$ is a $k$-algebra graded by an
abelian group $G$ (which we write additively), and that $c : G\times
G\rightarrow \kx$ is a 2-cocycle.

\definition{6.1} There exists a $G$-graded $k$-algebra $R'$ equipped with a
$G$-graded $k$-linear isomorphism
$$ \matrix R & \longrightarrow & R' \\
           r & \longmapsto & r' \endmatrix $$
such that, for all $\alpha,\beta \in G$, the multiplication
of arbitrary homogeneous elements $r' \in R'_\alpha$ and $s' \in
R'_\beta$ is given by the rule $r's'= c(\alpha,\beta)(rs)'$. Up to a
$G$-graded $k$-algebra isomorphism, $R'$ depends only on the
cohomology class of $c$ \cite{\ArScTa, p\. 888}. We call $R'$ the
{\it twist of $R$ by $c$\/} (cf\. \cite{\ArScTa, Section 3}), and we
refer to the above function $R \rightarrow R'$ as the {\it twist map
(associated to $c$)}. If $U$ is a subset of $R$, we will use $U'$ to
denote the image of $U$ under the twist map.

Note, for computational purposes, that if $I$ and $J$ are
($G$-)homogeneous ideals of $R$ then $I'J' = (IJ)'$.  \enddefinition

\definition{6.2} Recall that a proper homogeneous
ideal $K$, in $R$ or $R'$, is {\it graded-prime\/} provided it contains no
product $IJ$ of homogeneous ideals $I$ and $J$ not contained in
$K$. Every ideal of $R$ or $R'$ contains a unique maximum
homogeneous ideal, and an argument similar to (1.3) demonstrates that
the maximum homogeneous ideal within a prime ideal is graded-prime.
\enddefinition

\proclaim{Lemma} Suppose that $G$ is torsionfree and that $R$ is
noetherian. If $P_1,\dots,P_m$ are the minimal prime ideals of $R$,
then $P_1,\dots,P_m$ are homogeneous, and $P_1', \dots, P_m'$ are the
minimal prime ideals of $R'$. \endproclaim

\demo{Proof} Let $P_i^\circ$ denote the largest homogeneous ideal
contained in $P_i$. Since $G$ is torsionfree, the graded-prime ideal
$P_i^\circ$ must actually be a prime ideal, by \cite{\AbrHae,
Corollary 3.3}. Hence, $P_i= P_i^\circ$ by minimality.

Since the twist map preserves homogeneous ideals and their products,
it follows that $P_i'$ is a graded-prime ideal of $R'$. A second
application of \cite{\AbrHae, Corollary 3.3} then shows that $P_i'$ is
a prime ideal of $R'$. Now some product of the $P_i$ is equal to $0$
in $R$, and by homogeneity the corresponding product of the $P_i'$
must be $0$ in $R'$. It follows that every minimal prime ideal of $R'$
must occur among the $P_i'$.

Let $Q'_1,\dots,Q'_n$ be the minimal prime ideals of $R'$, and denote
their respective preimages under the twist map by $Q_1,\ldots,Q_n$.
Arguing as above, we see that the $Q_j$ are prime ideals of $R$, and
that all minimal prime ideals of $R$ occur among the $Q_j$. The lemma
follows. \qed\enddemo

\proclaim{6.3 Theorem} Let $k$ be an algebraically closed field, let
$R$ be a commutative affine $k$-algebra graded by a torsionfree
abelian group $G$, and let $A$ be the twist of $R$ by a 2-cocycle $c:
G\times G\rightarrow \kx$. Assume either that $-1$ is not in the
subgroup of $\kx$ generated by the image of $c$, or that $\chr
k=2$. Then there exist topological quotient maps $$\spec R\rightarrow
\spec A \qquad\text{and}\qquad \max R \rightarrow \prim A.$$
\endproclaim

\demo{Proof} Since $R$ is affine, its support is contained in a
finitely generated subgroup $G^{\text{fg}}\subseteq G$. Then $R$ is
also graded by $G^{\text{fg}}$, and $A$ is equal to the twist of $R$
by $c|_{G^{\text{fg}}}$. Hence, there is no loss of generality in
assuming that $G$ is finitely generated. Thus $G$ is now a free
abelian group of finite rank.

Set $\sigma(\alpha, \beta)= c(\alpha, \beta) c(\beta, \alpha)^{-1}$,
for $\alpha, \beta \in G$. By \cite{\ArScTa, Proposition 1(ii),
p\. 888}, $\sigma$ is an alternating bicharacter on $G$. Since either
$-1 \notin \langle \operatorname{im} \sigma \rangle$ or $\chr k=2$, it
follows from (4.2) that there exists an alternating bicharacter $d$ on
$G$ with the following properties: $d^2=\sigma$, and $d(\alpha,
\beta)=1$ whenever $\sigma( \alpha, \beta) =1$. By \cite{\ArScTa,
Proposition 1(i), p\. 888}, $d$ and $c$ are in the same cohomology
class, and hence $A$ can also be written as the twist of $R$ by
$d$. Thus, there is no loss of generality in assuming that $c$ is an
alternating bicharacter and that $c(\alpha, \beta) =1$ whenever
$c(\alpha, \beta)^2 =1$.

The inverse of the twist map is a $G$-graded $k$-linear isomorphism
$\Phi : A\rightarrow R$ such that $\Phi(ab)= c(\alpha, \beta) \Phi(a)
\Phi(b)$ for $a\in A_\alpha$ and $b\in A_\beta$. Choose homogeneous
elements $r_1,\dots,r_n$ generating $R$ as a $k$-algebra, and let
$\delta_1, \dots, \delta_n\in G$ be the degrees of these elements.

Set $\Gamma= \ZZ^n$ and $\Gamma^+= (\ZZ^+)^n$, and let $\rho : \Gamma
\rightarrow G$ be the group homomorphism given by
$$\rho(\alpha_1, \dots, \alpha_n)= \alpha_1 \delta_1 +\dots+ \alpha_n
\delta_n.$$
Set $\ctil(\alpha, \beta)= c(\rho(\alpha), \rho(\beta))$ for $\alpha,
\beta \in \Gamma$. Then $\ctil$ is an alternating bicharacter on $\Gamma$
such that $\ctil(\alpha, \beta) =1$ whenever $\ctil(\alpha, \beta)^2 =1$.
Now set $r^\alpha= r_1^{\alpha_1} r_2^{\alpha_2} \cdots r_n^{\alpha_n} \in
R_{\rho(\alpha)}$ and $a_\alpha= \Phi^{-1}(r^\alpha) \in A_{\rho(\alpha)}$
for $\alpha \in \Gamma^+$. Since $R$ is commutative, $r^\alpha r^\beta=
r^{\alpha +\beta}$ for $\alpha, \beta \in \Gamma^+$, whence
$$\Phi(a_\alpha a_\beta)= c(\rho(\alpha), \rho(\beta)) r^\alpha r^\beta=
\ctil(\alpha, \beta) r^{\alpha+ \beta}= \ctil(\alpha, \beta) \Phi(
a_{\alpha +\beta}), $$
and so $a_\alpha a_\beta= \ctil(\alpha, \beta) a_{\alpha +\beta}$ for
$\alpha, \beta \in \Gamma^+$.

Now set $\Rtil= \kgam^+$ and $\Atil= k^{\ctil}\Gamma^+$. Write $\Rtil$ in
terms of a $k$-basis $\{ y_\alpha \mid \alpha \in \Gamma^+ \}$ such that
$y_\alpha y_\beta= y_{\alpha +\beta}$ for $\alpha, \beta \in \Gamma^+$,
and write $\Atil$ in terms of a $k$-basis $\{ x_\alpha \mid \alpha\in
\Gamma^+ \}$ such that $x_\alpha x_\beta= \ctil(\alpha, \beta) x_{\alpha+
\beta}$ for $\alpha, \beta\in \Gamma^+$. Let $\phtil : \Atil \rightarrow
\Rtil$ be the $k$-linear isomorphism such that $\phtil(x_\alpha)=
y_\alpha$ for $\alpha \in \Gamma^+$. There are $k$-algebra epimorphisms $f
: \Atil \rightarrow A$ and $g : \Rtil \rightarrow R$ such that
$f(x_\alpha)= a_\alpha$ and $g(y_\alpha)= r^\alpha$ for $\alpha \in
\Gamma^+$. We now obtain the following commutative diagram:
$$\CD \Atil @>{\tilde{\Phi}}>> \Rtil\\
@V{f}VV @VV{g}V\\
A @>\Phi>> R \endCD$$
Thus if $I=\ker f$, then $\phtil(I)= \ker g$.

Define $\G$ as in (4.9), with $\ctil$, $\Rtil$, $\Atil$ playing the roles
of $c$, $R$, $A$. By (4.11), the rule $P\mapsto \phtil^{-1}(P:\G)$ gives
topological quotient maps
$$\phi_s : \spec\Rtil \rightarrow \spec\Atil \qquad\text{and}\qquad \phi_m
: \max\Rtil \rightarrow \prim\Atil.$$
These maps restrict to topological quotient maps
$$\gather \phi_s^{-1} V(I) \rightarrow V(I) \approx \spec A\\
\phi_m^{-1} \bigl( V(I)\cap \prim\Atil \bigr) \rightarrow V(I)\cap
\prim\Atil \approx \prim A. \endgather$$
To finish the proof of the theorem, it suffices to show that
$\phi_s^{-1}V(I)= V(\phtil(I))$. The inclusion $\subseteq$ is easy: If
$P\in \phi_s^{-1}V(I)$, then $\phtil^{-1}(P:\G) \supseteq I$, whence
$P\supseteq (P:\G)\supseteq \phtil(I)$.

Let $P_1,\dots,P_m$ be the prime ideals of $\Atil$ minimal over $I$. In
view of (6.2), the sets $\phtil(P_1), \dots, \phtil(P_m)$ are the prime
ideals of
$\Rtil$ minimal over $\phtil(I)$. By (4.10), each $\phtil(P_j) \in \ggspec
\Rtil$. Given any $P\in V(\phtil(I))$, we have $P\supseteq \phtil(P_j)$
for some $j$. Since $\phtil(P_j)$ is a $\G$-prime ideal, $(P:\G) \supseteq
\phtil(P_j) \supseteq \phtil(I)$, and thus $\phi_s(P)\in V(I)$. Therefore
$\phi_s^{-1}V(I)= V(\phtil(I))$, as desired. \qed\enddemo

\definition{6.4 Remark} Retain the notation in the statement of the
preceding theorem, and let $r_1,\ldots,r_n$ be homogeneous generators
for $R$. As seen in the proof, we may assume without loss of
generality that $G = \ZZ^n$ and that $c$ is an alternating bicharacter
on $\ZZ^n$. In particular, it can be shown that the topological
quotient maps in (6.3) are equivariant with respect to suitable
actions by a subgroup of $(\kx)^n$; details are left to the interested reader.
\enddefinition

\enddocument